\documentclass[11pt]{amsart}   
\usepackage{amsmath, amsthm, amssymb, amscd}                  

\DeclareMathOperator{\as}{asdim} 
\DeclareMathOperator{\ad}{ad} \DeclareMathOperator{\diam}{diam}
 \DeclareMathOperator{\dist}{dist}
 
\DeclareMathOperator{\Card}{Card}

\begin{document}

\newcommand{\To}{\longrightarrow}
\newcommand{\sA}{\mathcal{A}}
\newcommand{\sB}{\mathcal{B}}
\newcommand{\sU}{\mathcal{U}}
\newcommand{\sV}{\mathcal{V}}
\newcommand{\sW}{\mathcal{W}}
\newcommand{\sH}{\mathcal{H}}
\newcommand{\sY}{\mathcal{Y}}
\newcommand{\sZ}{\mathcal{Z}}
\newcommand{\Real}{\mathbf{R}}
\newcommand{\Complex}{\mathbf{C}}
\newcommand{\so}{\Rightarrow}
\newcommand{\sX}{\mathcal{X}}
\newcommand{\g}{\gamma}
\newcommand{\G}{\Gamma}
\newtheorem{Theorem}{Theorem}[section]
\newtheorem*{thm}{Theorem}
\newtheorem*{cor}{Corollary}
\newtheorem*{prop}{Proposition}
\newtheorem{Lemma}[Theorem]{Lemma}
\newtheorem{Proposition}[Theorem]{Proposition}
\newtheorem{Corollary}[Theorem]{Corollary}
\newtheorem*{defn}{Definition}
\newtheorem{Question}{Question}
\bibliographystyle{amsplain}

\title{Growth of the asymptotic dimension function for groups}
\author[G.C. Bell]{Gregory C. Bell}
\date{\today}
\address{Department of Mathematics, Penn State University,
University Park, PA 16802} \email{bell@math.psu.edu}


\subjclass[2000]{Primary: 20F69, Secondary 20E06, 20E22, 20F65}

\keywords{Asymptotic dimension, growth of dimension, graphs of
groups, relatively hyperbolic groups}

\begin{abstract}
It is relatively easy to construct a finitely generated group with
infinite asymptotic dimension: the restricted wreath product of
$\mathbb{Z}$ by $\mathbb{Z}$ provides an example.  In light of
this, it becomes interesting to consider the rate of growth of the
asymptotic dimension function of a group. Loosely speaking, we
measure the dimension on $\lambda$-scale and let $\lambda$
increase to infinity to recover the asymptotic dimension.  In this
paper we consider how the asymptotic dimension function is
affected by different constructions involving groups.
\end{abstract}
\maketitle
\section{Introduction}
The asymptotic dimension of a metric space was introduced by
Gromov \cite{Gr93} in his study of asymptotic invariants of
infinite groups.  Roughly speaking, the asymptotic dimension of a
metric space is the large-scale equivalent of covering dimension
of a topological space.  If $\Gamma$ is a finitely generated
group, one endows $\Gamma$ with a word metric associated to a
finite, symmetric generating set $S.$  Any choice $S'$ of finite,
symmetric generating set gives rise to a quasi-isometric metric
space (in fact, they are Lipschitz equivalent).  The asymptotic
dimension of a metric space is a quasi-isometry invariant, so we
can define the asymptotic dimension of a finitely generated group
$\Gamma$, $\as \Gamma$, as the asymptotic dimension of the metric
space corresponding to any finite, symmetric generating set.  More
generally, one can define $\as$ for coarse spaces; it turns out
that $\as$ is also a coarse invariant, see \cite{Ro03}.

Yu \cite{Yu1} showed that groups with finite asymptotic dimension
satisfy the Novikov higher signature conjecture.  Later, Yu
\cite{Yu2} generalized this result by proving that groups which
admit a uniform embedding into Hilbert space, in particular then,
groups with Yu's ``property A" satisfy the coarse Baum-Connes
conjecture, and hence the Novikov conjecture.

It is not difficult to construct a finitely generated group with
infinite asymptotic dimension.  Any group containing isomorphic
copies of $Z^n$ for each $n$ will have infinite asymptotic
dimension.  As pointed out by Roe in \cite{Ro03} two simple
examples are Thompson's group $F,$ described in \cite{Br96} and
the reduced wreath product of $\mathbb{Z}$ by $\mathbb{Z}.$

In this paper we consider the asymptotic dimension function of a
metric space $X.$  This function measures the dimension on the
scale $\lambda$ of the metric space $X.$  Although this function
is obviously not an invariant of quasi-isometry class of the
metric space, the growth of this function is, see Proposition
\ref{qi}.

Higson and Roe \cite{HR} showed that finitely generated groups
with bounded asymptotic dimension function have Yu's property A.
(By a theorem of Ozawa \cite{Oz}, Yu's property A for a finitely
generated group is equivalent to $C^\ast$-exactness of the group,
so this property is also often refered to as exactness of the
group.) Later, Dranishnikov \cite{Dr00} showed that bounded
geometry metric spaces whose asymptotic dimension function grows
to infinity sublinearly have property A. In \cite{Dr04},
Dranishnikov generalized this result, showing that groups whose
asymptotic dimension function grows at most polynomially have Yu's
property A, so the coarse Baum-Connes conjecture and Novikov
conjectures hold for such groups.

Here, we apply the techniques found in \cite{Be05,BD1,DG05,Os} and
others to show that the growth rate of the asymptotic dimension
function can be recovered from examining neighborhoods of
stabilizers of an isometric action of the group on a metric space
with finite asdim.  We apply this method to prove our main
theorems:\vskip10pt

\noindent\textbf{Corollary \ref{graphsofgroups}.}~\em Let $(G,Y)$
be a finite, connected graph of finitely generated groups.  Let
$\Gamma$ be the fundamental group of $(F,G).$  Then, the
asymptotic dimension function of $\Gamma$ grows no faster than the
asymptotic dimension function of the vertex groups $G_P$ in the
graph of groups.\rm\vskip10pt

Applying Osin's methods from \cite{Os} we are also able to
conclude:\vskip10pt

\noindent\textbf{Theorem \ref{relhyper}.}~\em Let $\Gamma$ be a
finitely generated group hyperbolic relative to a collection
$\{H_1,\ldots, H_n\}$ of subgroups. Then, the asymptotic dimension
function of $\Gamma$ grows no faster than the asymptotic dimension
function for each of the $H_i.$\rm \vskip10pt

We end the paper with some open questions about the asymptotic
dimension function.
\section{The asymptotic dimension function}

Let $(X,d)$ be a metric space.  Let $\sU$ be a cover of $X.$ A
\emph{Lebesgue number} for $\sU$ is a number $\lambda$ for which
every set $A\subset X$ with $\diam(A)\le \lambda$ is entirely
contained within a single element of $\sU.$  We denote the
Lebesgue number of $\sU$ by $L(\sU).$  The multiplicity of a cover
$\sU$ is $m(\sU)=\sup_{x\in X}\{\Card(\{U\in \sU\mid x\in U\})\}.$

We define the \emph{asymptotic dimension function} of the metric
space $X$ by
\[\ad_X(\lambda)=\min\{m(\sU)\mid L(\sU)\ge\lambda\}-1,\] where the
minimum is taken over all covers of $X$ by uniformly bounded sets.
Note that $\ad$ is monotonic and that
\[\lim_{\lambda\to\infty}\ad_X(\lambda)=\as X.\]

Let $f,g:\mathbb{R}_+\to\mathbb{R}_+.$  We write $f\preceq g$ if
there exists a $k\in\mathbb{N}$ so that $f(x)\le kg(kx+k)+k$ for
all $x\in\mathbb{R}_+.$  We write $f\approx g$ if $f\preceq g$ and
$g\preceq f.$

The metric spaces $(X_1,d_1)$ and $(X_2,d_2)$ are
\emph{quasi-isometric} if there exist constants $\lambda\ge
1,\epsilon\ge 0,$ and $C\ge 0$ and a map $f:X\to Y$ so that for
all $x$ and $y$ in $X_1,$
\[\frac1\lambda d_1(x,y)-\epsilon\le d_2(f(x),f(y))\le \lambda
d_1(x,y)+\epsilon,\] and every point of $X_2$ lies in the
$C$-neighborhood of the image of $f.$

It is easy to see that if $f:X_1\to X_2$ is a
$(\lambda,\epsilon)$-quasi-isometry then there is a
\emph{quasi-inverse} $f'$ for $f,$ i.e. a
$(\lambda',\epsilon')$-quasi-isometry $f':X_2\to X_1$ such that
there exists some $k$ for which $d(ff'(x'),x)\le k$ and
$d(f'f(x),x)\le k$ for all $x\in X_1$ and $x'\in X_2.$

We shall need the following simple construction to prove the
invariance of growth of $\ad$ for finitely generated groups.

\begin{Proposition}  Let $\sU$ be a cover of the metric space $X$
by uniformly bounded sets with multiplicity $m(\sU)$ and Lebesgue
number $L(\sU).$  Let $k\ll L(\sU).$  Then, there is a cover $\sV$
of $X$ by uniformly bounded sets with Lebesgue number $\ge
L(U)-2k$ and $k$-multiplicity $\le m(\sU).$
\end{Proposition}

\begin{proof}
Define $\sV$ to be the $-k$-neighborhood of $\sU,$ i.e.
\[V_\alpha=X\setminus(N_k(X\setminus U_\alpha)),\]
where $N_k(\cdot)$ denotes the $k$-neighborhood.  To see that
$\sV$ covers $X$ we take $x\in X.$  Since $k\ll L(\sU)$ there is a
$U\in \sU$ so that $B_{2k}(x)\subset U.$  Thus $x\in
X\setminus(N_k(X\setminus U)).$  The sets $V$ are certainly
uniformly bounded, so it remains only to check that
$k$-multiplicity of $V\le m(\sU).$  To this end, suppose $B_k(x)$
meets $V_1,\ldots,V_p.$  Since $N_k(V_i)\subset U_i$ for each $i,$
we see that $x\in N_k(V_i)$ implies that $x\in U_i.$  Thus, $p\le
m(\sU).$  Finally, observe that $\diam(A)<L(\sU)-2k$ implies that
$\diam(N_k(A))<L(\sU).$  So $N_k(A)\subset U$ for some $U\in\sU.$
Thus $A\subset X\setminus(N_k(X\setminus U)).$
\end{proof}

Dranishnikov remarks in \cite{Dr04} that the growth of $\ad$ is a
quasi-isometry invariant and that this was known to Gromov. In
this paper we are only interested in the growth of the asymptotic
dimension function for groups, so we prove a weaker version,
showing the growth is an invariant of quasi-isometry when the
spaces are discrete with bounded geometry.

\begin{Proposition} \label{qi} Let $X$ and $Y$ be discrete metric spaces with
bounded geometry. Suppose that $X$ and $Y$ are quasi-isometric.
Then $\ad_X\approx\ad_Y.$ In particular, the $\approx$-equivalence
class of $\ad_\Gamma$ is well-defined for finitely generated group
$\Gamma.$
\end{Proposition}

\begin{proof}  Let $f:X\to Y$ be an
$(\alpha,\epsilon)$-quasi-isometry with
$(\alpha',\epsilon')$-quasi-inverse $g:Y\to X.$  We have to show
that there is some $k>0$ so that $\ad_Y(\lambda)\le
k\ad_X(\lambda+k)+k.$  Take $C$ so large that $N_C(f(X))\supset Y$
and both $d(fg(y),y)< C$ and $d(gf(x),x)< C$ for all $x\in X$ and
all $y\in Y.$

Let $\lambda\gg C$ be given.  Suppose that $\sU$ is a cover of $X$
by uniformly bounded sets such that $\sU$ has Lebesgue number
$L(\sU)>\alpha'\lambda+(\epsilon'+2C)$ and multiplicity
$1+\ad_X(\alpha'\lambda+(\epsilon'+2C)).$  Use the previous
proposition to define a cover $\bar\sU$ of $X$ with Lebesgue
number $>\alpha'\lambda+\epsilon'$ and $C$-multiplicity
$\le\ad_X(\alpha'\lambda+(\epsilon'+2C))+1.$  Put
$\sV=\{N_C(f(U))\mid U\in\bar\sU\}.$  The collection $\sV$ forms a
cover of $Y$ by the choice of $C.$  It is easy to see that $\sV$
consists of uniformly bounded sets, so it remains to compute its
Lebesgue number and its multiplicity.

Suppose that $B\subset Y$ with $\diam(B)\le\lambda.$  Then
$\diam(g(B))\le\alpha'\lambda+\epsilon'.$  So there is some $U\in
\bar \sU$ so that $g(B)\subset U$ and hence $fg(B)\subset f(U).$
If $b\in B,$ then $d(b,fg(b))< C,$ so $b\in N_C(U)$ which is an
element of $\sV.$  Thus, $L(\sV)>\lambda.$

Finally, suppose $y\in Y$ meets $V_1,\ldots,V_p$ in $\sV.$ Then
$B_C(y)$ meets each $f(U_i),$ say $d(y,f(x_i))<C$ for all
$i=1,2,\ldots, p$ with $x_i\in U_i.$  We claim that $p\le
c_Y(C)(\ad_X(\alpha'\lambda+\epsilon'+2C)+1),$ where $c_Y(C)$ is
the constant from the bounded geometry condition on $Y.$  First
observe that $B_C(y)$ contains at most $c_Y(C)$ distinct points.
If one such point, say $y_0$ were the $f$-image of
$x_0^1,\ldots,x_0^t$ where $x_0^i$ are in distinct members $U_i$
of the collection $\bar\sU,$ then we see that $d(g(y_0),x_0^i)\le
C$ for all $i=1,2,\ldots,t.$  Since the $C$-multiplicity of
$\bar\sU$ is bounded above by
$\ad_X(\alpha'\lambda+\epsilon'+2C)+1,$ we see that $p\le
c_Y(C)(\ad_Y(\alpha'\lambda+\epsilon'+2C)+1)$ as desired.

Since $\ad$ is monotonic, by setting
$k>\max\{c_Y(C),\alpha',\epsilon'+2C\}$ we see that
$\ad_Y(\lambda)\le k\ad_X(k\lambda+k)+k,$ as required.
\end{proof}

Often we will apply the preceding result to an $R$-neighborhood
$N_R(A)$ of a set $A$ in a finitely generated group $\Gamma$ to
conclude that $\ad_{N_R(A)}\approx \ad_A.$

\section{Groups acting on finite-dimensional spaces}

We will say that a family $\{X_\alpha\}_{\alpha}$ has
\emph{uniform asymptotic dimension growth} if for every $\lambda$
one can find a constant $B=B(\lambda)$ and a family
$\{\sU_\alpha^\lambda\}$ of covers, $\sU_\alpha$ covering
$X_\alpha$ so that $B$ forms a uniform bound on the diameters of
the sets in each cover, so that the Lebesgue number of
$\sU_\alpha$ exceeds $\lambda$ for all $\alpha,$ and so that the
multiplicity of each family is uniformly bounded in $\alpha.$  It
is clear that any finite collection of subsets of a metric space
satisfies these conditions, as does a collection of subsets
belonging to the same isometry class.

\begin{Theorem} \label{union} Let $X=\bigcup_\alpha X_\alpha$ be a metric space such
that the collection $\{X_\alpha\}$ has uniform asymptotic
dimension growth $\ad(\lambda).$  Moreover, for each $\lambda,$
assume that there is some set $Y_\lambda\subset X$ so that
$\{X_\alpha\setminus Y_\lambda\}$ forms a $3B(\lambda)$-disjoint
collection and that $Y_\lambda$ is covered by a family of
uniformly bounded sets $\sV$ with $L(\sV)>\lambda$ and
multiplicity $\le m_Y(\lambda).$  Then $ad_X(\lambda)\approx
\max\{ad(\lambda),m_Y(\lambda)\}.$
\end{Theorem}

\begin{proof}  Let $\lambda>0$ be given.
For each $\alpha,$ let $\sU^\lambda_\alpha$ be a cover of
$X_\alpha$ so that the collection $\{\sU^\lambda_\alpha\}_\alpha$
satisfies the the uniform condition for asymptotic dimension
growth, with $B(\lambda)\gg\lambda.$ Define
$\bar\sU_\alpha^\lambda$ to be the collection $\{U\setminus
N_{-\lambda}(Y_\lambda)\mid U\in \sU^\lambda_\alpha\},$ and put
$\sW=\sW^\lambda =\bigcup_{\alpha}\bar\sU^\lambda_\alpha \cup
\sV.$ We claim that $\sW$ covers $X,$ is uniformly bounded,
$L(\sW)>\lambda$ and that $m(\sW)\le \ad(\lambda)+m_Y(\lambda).$

It is obvious that $\sW$ covers $X.$

Next, we have $\diam(\bar\sU^\lambda_\alpha)\le B(\lambda).$ Since
$\diam(\sW)\le\diam(\bar\sU^\lambda_\alpha)+\diam(\sV),$ the cover
$\sW$ consists of uniformly bounded sets.

Next, let $x\in X.$  First, we claim that there can be at most one
index $\beta$ with $x$ in sets from $\bar\sU^\lambda_{\beta}.$
Suppose that $x\in U$ and $x\in U'$ with
$U\in\bar\sU^\lambda_\alpha$ and $U'\in\bar\sU^\lambda_{\alpha'}.$
Then $x\in X_\alpha\setminus N_{-\lambda}(Y_\lambda)\cap
X_{\alpha'}\setminus N_{\lambda}(Y_{\lambda}),$ but
$\dist(X_\alpha\setminus
N_{-\lambda}(Y_\lambda),X_{\alpha'}\setminus
N_{-\lambda}(Y_\lambda))\ge 3B(\lambda)-2\lambda\ge\lambda.$ Thus,
there can be at most $m_Y(\lambda)$ sets from $\sV$ containing $x$
and at most $\ad(\lambda)$ sets of the form
$\bar\sU^\lambda_\alpha$ containing $x.$  So $m(\sW)\le
m_Y(\lambda)+m(\sU^\lambda_\alpha).$

Finally, we check that the Lebesgue number is $>\lambda.$  To this
end, let $A\subset X$ with $\diam(A)\le\lambda.$  There are only
two possibilities for $A.$  Either it lies entirely within
$Y_\lambda,$ in which case, it is contained in a set from $\sV,$
or it is contained in $X_\alpha\setminus
N_{-\lambda}(Y_{\lambda}).$  In the second case it is contained in
one of the elements of $\bar\sU^\lambda_\alpha.$
\end{proof}

Suppose $X=A\cup B.$ Setting $Y_\lambda=B$ the disjointness
condition is trivially satisfied so we immediately obtain:

\begin{Corollary} \emph{(Finite Union Theorem)} Let $X=A\cup B$ be a metric space with
asymptotic dimension functions $ad_A(\lambda)$ and
$ad_B(\lambda),$ respectively.  Then
$ad_X(\lambda)\approx\max\{ad_A(\lambda),ad_B(\lambda)\}.$
\end{Corollary}

A version of the union theorem we have occasion to use frequently
is the following:

\begin{Corollary} \label{unioncor} Let $X$ be a metric space with bounded geometry
(or a finitely generated group with a word metric). Suppose
$X_\alpha$ is a collection of isometric subsets of $X$ and that
$\bigcup_\alpha X_\alpha=X.$ Suppose for each $r>0,$ there exists
a $Y_r\subset X$ with $\{X_\alpha\setminus Y_r\}$ $r$-disjoint,
where $Y_r$ is quasi-isometric to $X_\alpha.$ Then
$\ad_{X}\approx\ad_{X_\alpha}.$
\end{Corollary}

\begin{proof} Since the $X_\alpha$ are isometric, their asymptotic
dimension functions grow uniformly as $\ad_{X_\alpha}.$ Let
$\lambda>0$ be given and take $r\ge 3B(\lambda),$ where
$B(\lambda)$ is the constant from the definition of uniform
asymptotic dimension function growth.  Then $\{X_\alpha\setminus
Y_r\}$ is $3B(\lambda)$ disjoint.  Since $Y_r$ is quasi-isometric
to any element $X_\alpha$ in the collection, there is a uniformly
bounded cover of $Y_r$ with Lebesgue number $>\lambda$ and with
multiplicity $m_Y(\lambda)\le k\ad_{X_\alpha}(k\lambda+k)+k$ for
some $k\in\mathbb{N}.$  So, by the union theorem,
$ad_X\approx\ad_{X_\alpha}.$
\end{proof}

Observe that the same result holds if the set $Y_r$ is assumed to
be a finite union of spaces each of which is quasi-isometric to an
$X_\alpha.$

\begin{Theorem} \label{main} Let $X$ be a metric space with finite
asymptotic dimension. Suppose that the finitely generated group
$\Gamma$ acts by isometries on $X.$ Finally, suppose there is some
$f:\mathbb{R}_+\to\mathbb{R}_+$ so that for every $R>0,$
$ad_{W_R}\preceq f.$  Then, $ad_{\Gamma}\preceq f.$
\end{Theorem}

\begin{proof} Let $S=S^{-1}$ be a finite generating set for
$\Gamma.$  Fix a point $x_0\in X,$ and define $\pi:\Gamma\to X$
via the action $\gamma\mapsto\gamma.x_0.$  Let
$\mu=\max\{\dist_X(s.x_0,x_0)\mid s\in S\}.$  We claim that $\pi$
is $\mu$-Lipschitz.  Since the metric on $\Gamma$ is discrete
geodesic, it suffices to verify the Lipschitz condition on pairs
$(\gamma,\gamma'),$ where $\dist_S(\gamma,\gamma')=1.$  Such a
pair is necessarily of the form $(\gamma,\gamma s),$ where $s\in
S.$ Computing, we find $\dist_X(\pi(\gamma),\pi(\gamma
s))=\dist_X(\gamma.x_0,\gamma s.x_0)=\dist_X(x_0,s.x_0)\le\mu.$

Let $R>0$ and $\lambda>0$ be given.  By assumption $\as X<\infty,$
so put $\as X=k-1.$ There exists an $R/2$-uniformly bounded cover
$\sU$ of $\Gamma.x_0$ with Lebesgue number $\ge\lambda\mu$ and
multiplicity $\le k.$

By assumption there is a uniformly bounded cover $\sV$ of
$W_R(x_0)$ with Lebesgue number $\ge \lambda$ and multiplicity
$\le nf(n\lambda+n)+n,$ for some $n.$  Denote $nf(n\lambda+n)+n$
by $f_n(\lambda).$ For each $U\in\sU,$ let $\gamma_U\in\Gamma$ be
an element such that $\gamma_U x_0\in\sU.$ Since
left-multiplication is an isometry on $\Gamma,$ we see that
$\gamma_U W_R(x_0)$ is isometric to $W_R(x_0).$  Thus, we can take
covers of the form $\gamma_UV$ for $\gamma_U W_R(x_0)$ that are
uniformly bounded, have multiplicity $\le f_n(\lambda)$ and have
Lebesgue number $\ge \lambda.$

Define a family $\sW$ of subsets of $\Gamma$ by
\[\sW=\{\gamma_UV\cap\pi^{-1}(U)\mid U\in\sU,V\in\sV\}.\]
We claim that $\sW$ is a uniformly bounded cover of $\Gamma$ with
Lebesgue number $\ge \lambda$ and multiplicity $\le
kf_n(\lambda).$

If $\gamma\in\Gamma$ then there is some $U\in\sU$ so that
$\gamma.x_0\in U.$  Thus, $\dist_X(\gamma_U.x_0,\gamma.x_0)\le R.$
So $\gamma\in\pi^{-1}(B_R(\gamma_U.x_0))$ and we conclude that
there is some $V\in\sV$ for which $\gamma_UV$ contains $\gamma.$
Thus, $\gamma\in\pi^{-1}(U)$ and $\gamma\in\gamma_UV,$ so there is
a $W\in\sW$ containing $\gamma.$

Next, let $A\subset\Gamma$ be given with $\diam_\Gamma(A)\le
\lambda.$ Then, we have $\diam_X(\pi(A))\le\mu\lambda.$  Since
$L(\sU)\ge \lambda\mu,$ there is a $U\in\sU$ so that
$\pi^{-1}(U)\supset A.$ Next,
$\diam_\Gamma(\gamma_U^{-1}(A))\le\lambda,$ and
$\gamma_U^{-1}(A)\subset W_R(x_0),$ so there is a $V\in\sV$ with
$\gamma^{-1}_U(A)\subset V.$  Thus, $\gamma_UV$ will contain $A,$
and so there is an element of $\sW$ containing $A.$

Let $\gamma\in\Gamma.$  Then $\pi(\gamma)$ is in at most $k$ of
the sets in $\sU.$  On the other hand, in each $\gamma_U\sV,$
$\gamma$ can belong to at most $f_n(\lambda)$ of the $\gamma_UV.$
Thus, the multiplicity of $\sW$ is bounded above by
$kf_n(\lambda).$

The family $\sW$ is uniformly bounded since $\sV$ is and since
left multiplication is an isometry in $\Gamma.$
\end{proof}

Using the techniques found in \cite{BD1} it is easy to prove the
following result, (compare \cite[Section 4]{Dr04}):

\begin{Corollary} Let $\phi:G\to H$ be a surjective homomorphism
of finitely generated groups with $\ker\phi=K.$  If $\as
H<\infty,$ then $\ad_G\approx\ad_K.$
\end{Corollary}

As an extension theorem for the asymptotic dimension function of
groups, this is unsatisfactory due to the requirement that the
group $H$ must have finite $\as.$  A much more satisfactory
statement would be: $\ad_G\approx\max\{\ad_H,\ad_K\},$ but our
techniques do not apply to this situation.

\section{Graphs of groups and complexes of groups}

The Bass-Serre theory of graphs of groups Bass-Serre \cite{Se} is
well known. If $Y$ is a finite, connected graph we can label the
vertices of $Y$ with groups.  An edge joining two vertices should
be equipped with a group as well as two injective homomorphisms
from the edge group into the two vertex groups connected by the
edge.  One then forms a group (the \emph{fundamental group of the
graph of groups}) by taking the free product of the vertex groups
and the edge labels (with formal inverses) and requiring certain
compatibility conditions on the group based on the inclusions of
the edge groups into the vertex groups.  For more details, see
below or \cite{BH,Se}.

The group constructed in this way will then act on a tree by
isometries and the quotient of the tree by this action is the
original graph $Y.$  A nice feature of the theory is that to every
graph of groups there is an associated group and group action, and
for every action by isometries of a group on a tree one can
recover the original group via the graph of groups construction.

Three standard examples of fundamental groups of graphs of groups
are the following:
\begin{itemize}
    \item If $Y$ is any graph, and all vertex groups are taken to
    be the trivial group, then the fundamental group of the graph
    of groups is the fundamental group $\pi_1(Y).$
    \item If $Y$ has two vertices and one edge, then the
    fundamental group of the graph of groups is the free product
    of the vertex groups amalgamated over the edge group.
    \item If $Y$ has one vertex and one edge, then the fundamental
    group of the graph of groups is the HNN-extension of the
    single vertex group over the edge group.
\end{itemize}

There is a natural generalization of the Bass-Serre theory of
graphs of groups (called complexes of groups) due to Haefliger
\cite{Hae}. An exposition of this theory can be found in
\cite{BH}; we follow the notation found there.  Roughly speaking,
one replaces the action on a tree with an action on a suitable
higher-dimensional replacement called a \textbf{s}mall
\textbf{c}ategory \textbf{w}ith\textbf{o}ut \textbf{l}oops (called
a scwol).  One issue that arises with the theory of complexes of
groups is that there is not always an associated isometric action
of the fundamental group of a complex of groups on scwol so that
the quotient is the original complex.  When this happens, the
complex of groups is said to be \emph{developable}.  Since we are
interested in exploiting this action to our advantage, we only
consider developable complexes of groups.  In the language of
complexes of groups, the Bass-Serre theory states that complexes
of groups associated to $1$-dimensional scwols (see below) are
developable.

The notation we need and pertinent results from the theory of
complexes of groups follow.

\begin{defn}  A {\em small category without loops} (abbreviated
{\em scwol}) is a set $\sX$ which is the disjoint union of a {\em
vertex set} $V(\sX)$ and an {\em edge set} $E(\sX).$  There are
maps \[i:E(\sX)\to V(\sX)\qquad\text{and}\qquad t:E(\sX)\to
V(\sX)\] which assign to each edge $a$ the initial vertex of $a$
and the terminal vertex of $a,$ respectively.  Let $E^{(k)}(\sX)$
denote the composable sequences of edges of length $k,$ i.e.,
$E^{(k)}(\sX)=\{(a_1,\ldots,a_k)\in (E(\sX))^k\mid
i(a_i)=t(a_{i+1}),$ for $i=1,\ldots, k-1\}.$  By convention,
$E^{(0)}(\sX)=V(\sX).$ There is also a map
\[E^{(2)}(\sX)\to E(\sX)\] which assigns to each pair $(a,b)$ an
edge $(ab)$ called the composition of $a$ and $b.$  These maps are
required to satisfy:
\begin{enumerate}
\item $i(ab)=i(b),$ and $t(ab)=t(a)$ for all $(a,b)\in
E^{(2)}(\sX);$ \item $a(bc)=(ab)c$ for all edges $a,b$ and $c$
with $i(a)=t(b)$ and $i(b)=t(c);$ and \item $i(a)\neq t(a).$ (the
no loops condition)
\end{enumerate}
\end{defn}

We define the {\em dimension} of the scwol $\sX$ to be the maximum
$k$ such that $E^{(k)}(\sX)$ is not empty.

\begin{defn}  The {\em geometric realization} $|\sX|$ is a piecewise
Euclidean polyhedral complex, with each $k$-cell isometric to the
standard simplex $\Delta^k.$  There is one such $k$-simplex $A$
for each $A\in E^{(k)}(\sX).$  The identifications are the obvious
ones, induced by the face relation among simplices.
\end{defn}

The geometric realization is a Euclidean complex and
is given its intrinsic metric.  

\begin{defn}  A {\em group action} on a scwol is a homomorphism $G\to
Aut(\sX)$ satisfying
\begin{enumerate}
\item for every $g\in G,$ and for all $a\in E(\sX)$ $g.i(a)\neq
t(a).$ \item for every $g\in G,$ and for all $a\in E(\sX)$ if $g.
i(a)= i(a),$ then $g. a=a.$
\end{enumerate}
\end{defn}

Notice that a group action on a scwol induces an isometric action
of the group on the geometric realization $|\sX|.$  Since we are
primarily concerned with isometric actions on metric spaces, this
is the action that we consider.

One forms the quotient $\sY=G\backslash \sX$ of the scwol $\sX$ by
the action of $G$ by taking $V(\sY)=G\backslash V(\sX),$ and
$E(\sY)= G\backslash E(\sX).$  One can verify that $\sY$ has the
structure of a scwol.

\begin{defn}  A {\em complex of groups $G(\sY)$ over a scwol} $\sY$ is a
collection $G(\sY)=(G_\sigma, \psi_a, g_{a,b})$ satisfying
\begin{enumerate}
\item to each $\sigma\in V(\sY),$ there corresponds a group
$G_\sigma$ called the {\em local group} at $\sigma;$ \item  for
each $a\in E(\sY)$ there exists an injective homomorphism
$\psi_a:G_{i(a)}\to G_{t(a)};$ and \item For each $(a,b)\in
E^{(2)}(\sY),$ there is a $g_{a,b}\in G_{t(a)}$ such that
\newline\indent\indent(i) $Ad(g_{a,b})\psi_{ab}=\psi_a\psi_b,$
where $Ad(g_{a,b})$ denotes conjugation by $g_{a,b},$ and
\newline\indent\indent (ii)
$\psi_a(g_{b,c})g_{a,bc}=g_{a,b}g_{ab,c},$ for all $(a,b,c)\in
E^{(3)}(\sY).$
\end{enumerate}
\end{defn}


When a complex of groups is developable, there is an explicit
method of constructing both the scwol $\sX$ and the group $G$
which acts on the scwol.  The scwol $\sX$ on which the group acts
is simply connected and has an explicit description in a similar
way to the construction of the tree $\tilde X$ in the theory of
graphs of groups (see \cite{Se}).

Indeed, if $G(\sY)$ is a developable complex of groups, then we
can define the development $D(\sY)$ to be the scwol whose vertices
and edges are given by $V(D(\sY))=\{(gG_\sigma, \sigma)\mid
\sigma\in V(\sY)\},$ and $E(D(\sY))=\{(gG_{i(a)},a)\mid a\in
E(\sY)\}.$ Then the group $G$ acts on the development $D(\sY)$ by
left multiplication. The development is isomorphic to the scwol
$\sX,$ mentioned above.  (See \cite{BH} for more details.)

We describe the fundamental group of the complex of groups
$\pi_1(G(\sY))$ which is the group $G,$ up to isomorphism.  As in
the theory of graphs of groups, there are two equivalent
descriptions of the fundamental group, but we only describe one.
It relies on the construction of the auxiliary group $FG(\sY).$
Let $E^\pm(\sY)$ denote the collection of symbols $\{a^+,a^-\}$
where $a\in E(\sY).$  The elements of $E^\pm(\sY)$ can be thought
of as {\em oriented edges}.  If $e=a^+,$ then define $i(e)=t(a)$
and $t(e)=i(a).$  Accordingly, if $e=a^-,$ define $t(e)=t(a)$ and
$i(e)=i(a).$ Then, define $FG(\sY)$ to be the free product of the
local groups $G_\sigma$ and the free group generated by the
collection $E^\pm(\sY)$ subject to the additional relations:
\begin{enumerate}
\item $(a^+)^{-1}=a^-,$ and $(a^-)^{-1}=a^+;$ \item
$a^+b^+=g_{a,b}(ab)^+;$ \item $\psi_a(g)=a^+ga^-,$ for all $g\in
G_{i(a)}.$
\end{enumerate}

An edge path in $\sY$ is
a sequence $(e_1,\ldots,e_k)$ with $t(e_i)=i(e_{i+1}),$ for all
$i=1,\ldots,k-1.$ By a $G(\sY)$-path issuing from the vertex
$\sigma_0$ we mean a sequence $(g_0,e_1,g_1,\ldots,e_k,g_k),$
where $(e_1,\ldots,e_k)$ is an edge path in $\sY,$ $g_0\in
G_{\sigma_0},$ $i(e_1)=\sigma_0,$ and $g_i\in G_{t(e_i)},$ for
$i>0.$  We associate the word $g_0e_1\ldots e_kg_k\in FG(\sY)$ to
the path described above. A $G(\sY)$-loop based at $\sigma_0$ is a
$G(\sY)$ path with $t(e_k)=\sigma_0.$ The $G(\sY)$-loops $\gamma$
and $\gamma'$ are homotopic if the $FG(\sY)$-words they represent
are equal. The fundamental group $\pi_1(G(\sY),\sigma_0)$ is the
collection of all words associated to $G(\sY)$-loops based at
$\sigma_0,$ up to homotopy equivalence.


Let $G(\sY)$ be a developable complex of groups.  Fix a vertex
$\sigma_0$ in $\sY$ and consider the action of the fundamental
group $\pi=\pi_1(G(\sY),\sigma_0)$ on the simply connected scwol
$\sX_e$ induced by the complex of groups. In \cite[Proposition
1]{Be05} the $R$-stabilizers of this action are characterized as
being those elements with associated path $c$ of length not
exceeding $R.$

\begin{Lemma} Let $\pi$ denote the fundamental group of a complex
of groups $G(\sY)$ where $\sY$ is finite and connected, and the
local groups are finitely generated.  Let $\sigma_0$ be a vertex.
Suppose there is some function $f$ so that each stabilizer
$\Gamma_{\sigma}$ satisfies $ad_{\Gamma_\sigma}(\lambda)\le
f(\lambda).$ Then $ad_{W_R}\preceq f$ for all $R.$
\end{Lemma}

\begin{proof}  The proof of this result is analogous to
that of \cite[Lemma 2]{Be05}.

Let $K\subset FG(\sY)$ denote the set of all words in $FG(\sY)$
with an associated path issuing from $\sigma_0.$  Observe that
$\pi\subset K$ and the set $K$ acts on $\sX$ by left
multiplication.  We consider $W_R(\sigma_0)$ as a subset of $K$
and show that it has $\ad_{W_R}\preceq f.$  It will follow that as
a subset of $\Gamma,$ $\ad_{W_R}\preceq f.$

Put $K_j$ equal to the subset of $K$ whose associated paths have
length exactly equal to $j.$  Then, in light of the finite union
theorem, we need only show that $\ad_{K_j}\preceq f$ for all $j.$
We proceed inductively.  Observe that $K_0=G_{\sigma_0},$ so this
is true by assumption. Next, consider the case $K_{j+1}$ with
$j\ge 0.$ Observe that $K_{j+1}\subset \bigcup_{a\in E^\pm(\sY)}
K_j aG_{t(a)}.$

The orientation of the edge $a$ is an issue since it determines
whether the group $G_{t(a)}$ is a domain or codomain of the
function $\psi_a.$  Thus, it is necessary to consider two cases
separately.

Suppose first that $a$ has negative orientation.  So, we are
considering $K_ja^-G_{t(a)}.$  For every $r>0$ let
$Y_r=K_ja^-N_r(\psi_a(G_{i(a)})),$ where the $r$-neighborhood is
taken in the group $FG(\sY).$  Then $Y_r$ is quasi-isometric to
$K_ja^-\psi_a(G_{i(a)}).$  We observe that
$K_ja^-\psi_a(G_{i(a)})=K_jG_{i(a)} a^-,$ which is just $K_ja^-.$
Since $K_ja^-$ is quasi-isometric to $K_j,$ we have
$\ad_{Y_r}\approx\ad_{K_j},$ which by the inductive hypothesis
grows no faster than $f.$

Next, we decompose the set $K_ja^-G_{t(a)}$ into sets of the form
$\{xa^-G_{t(a)}\},$ where the index runs over all $x\in K_j$ that
do not end with an element $g\in G_{i(a)}.$  One can still obtain
these elements through the relations of $FG(\sY).$  For example,
to obtain $xga^-g^\prime,$ with $x$ of the required form, $g\in
G_{i(a)}$ and $g^\prime$ in $G_{t(a)},$ simply take the word
$xa^-\psi_a(g)g^\prime,$ which is of the required form. Next,
observe that the map $G_{t(a)}\mapsto xa^-G_{t(a)}$ is an isometry
in the (left-invariant) word metric.  So the set
$\{xa^-G_{t(a)}\}$ has $\ad \preceq f,$ uniformly.

In order to apply the union theorem to this family, it remains to
show only that the family $\{xa^-G_{t(a)}\setminus Y_r\}$ is
$r$-disjoint.  To this end, let $xa^-z$ and $x'a^-z'$ be given in
different families.  Then we compute
$d(xa^-z,x'a^-z')=\|z^{-1}a^+x^{-1}x'a^-z'\|.$ Since $z$ and $z'$
lie outside of $N_r(\psi_a(G_{i(a)})),$ take $z=\psi_a(g)s$ and
$z'=\psi_a(g')s',$ where $\|s\|>r,$ $\|s'\|>r,$ and $\psi_a(g)$
and $\psi_a(g')$ are in $\psi_a(G_{t(a)}).$  Then,
\[\|z^{-1}a^+x^{-1}x'a^-z'\|=\|s^{-1}a^+g^{-1}x^{-1}x'g'a^-s'\|.\]
Now, in order for this length to be less than $r,$ a reduction
must occur in the middle, so that $a^+$ and $a^-$ annihilate each
other.  In order for this to occur, we must have
$g^{-1}x^{-1}x'g'\in G_{i(a)}.$  Thus, $x^{-1}x'\in G_{i(a)}.$
But, this means that $xa^-G_{t(a)}$ and $x'a^-G_{t(a)}$ define the
same set.  Thus, in the case that the edge has negative
orientation, we have $\ad_{K_ja^-G_{t(a)}}\preceq f.$

Next, we consider the case where the edge $a$ has positive
orientation.  In this case,
$K_ja^+G_{i(a)}=K_j\psi_a(G_{i(a)})a^+,$ which is quasi-isometric
to $K_j.$ We conclude that
$\ad_{K_ja^+G_{i(a)}}\preceq\ad_{K_j}\preceq f.$
\end{proof}

It is clear that scwols of dimension 1 must have precisely two
types of vertices: sources and sinks.  A source is an initial
vertex of every edge it is contained in and a sink is a terminal
vertex of every edge it is contained in. Every one-dimensional
simplicial complex (graph) can be given the structure of a
one-dimensional scwol by placing a source vertex in the middle of
every edge, thus giving the original vertices the structure of
sinks.  It is easy to verify that the theory of complexes of
groups over one-dimensional scwols is precisely the same as the
theory of graphs of groups. Phrased in terms of the language of
complexes of groups, the Bass-Serre structure theorem for groups
acting without inversion on graphs says that if $dim(\sY)=1,$ then
$G(\sY)$ is always developable.

Phrasing this result in the language of graphs of groups, we
obtain:
\begin{Corollary} Suppose the finitely generated group $\Gamma$ acts by isometries on a tree
$X$ with compact quotient and finitely generated stabilizers such
that there is some function $f$ so that each stabilizer $\Gamma_x$
satisfies $ad_{\Gamma_x}(\lambda)\le f(\lambda).$ Then
$ad_{W_R}\preceq f$ for all $R.$
\end{Corollary}

Applying Theorem \ref{main} we obtain our main result for
developable complexes of groups:

\begin{Theorem} Let $\pi$ be the fundamental group of a finite,
developable complex of groups $G(\sY)$ such that the development
$\sX$ has $\as|\sX|<\infty,$ and such that there is some function
$f$ that is an upper bound for the growth of the asymptotic
dimension function of every base group $G_\sigma.$  Then
$ad_\Gamma\preceq f.$
\end{Theorem}

In the language of graphs of groups this becomes:

\begin{Theorem} \label{graphsofgroups} Let $\pi$ be the fundamental group of a finite
graph of groups.  Suppose that the vertex groups are finitely
generated and that the growth of the asymptotic dimension of the
vertex groups is bounded above by the function $f.$  Then,
$\ad_\pi\preceq f.$
\end{Theorem}

Based on the examples cited above we immediately obtain:
\begin{Corollary} Let $A$ and $B$ be finitely generated groups.
Then we have the following estimates on growth of $\as$ for the
amalgamated free product and HNN-extension: $\ad_{A\ast_C
B}\approx\max\{\ad_A,\ad_B\}$ and $\ad_{A\ast_C}\approx\ad_A.$
\end{Corollary}

\section{Relatively hyperbolic groups}

In a recent article, Osin \cite{Os} proved that a finitely
generated group that is hyperbolic relative to a finite collection
$\{H_1,\ldots, H_m\}$ of subgroups has finite asymptotic dimension
if $\as H_i<\infty$ for each subgroup $H_i.$  We refer the reader
to Osin's articles \cite{Os,Os1} among others for definitions of
relatively hyperbolic groups.

Dadarlat and Guentner \cite{DG05} use similar techniques to prove
that the group $\Gamma$ is uniformly embeddable in Hilbert space
precisely when each subgroup $H_i$ is uniformly embeddable in
Hilbert space. Using altogether different techniques, Ozawa
\cite{Oz05} was able to show the corresponding result for
$C^\ast$-exactness of such a group.  Ozawa's result is recovered
in the work of Dadarlat-Guentner.

Our goal in this section is to prove the corresponding result for
the growth of the asymptotic dimension function for the group
$\Gamma.$ Osin gives most of the ingredients and techniques for
this result, we simply put them together with our results on the
growth of the asymptotic dimension function.

First we fix some notation.  Throughout this section $\Gamma$ will
be a finitely generated group with generating set $S=S^{-1}$ that
is hyperbolic with respect to the collection $H_1,\ldots,H_m$ of
subgroups. Put $\sH=\cup_{i=1}^m(H_i\setminus\{e\}).$ There are
two metrics we wish to consider on $\Gamma.$  The first is the
word metric $d_S$ associated to the generating set $S.$  The
second will be denoted $d_{S\cup\sH}$ and it is the
(left-invariant) word metric associated to the set $S\cup\sH.$

Observe that there is an action of $\Gamma$ on
$(\Gamma,d_{S\cup\sH})$ by isometries.  Osin shows that
$\as(\Gamma,d_{S\cup\sH})<\infty.$  So, in order to apply Theorem
\ref{main} it remains to examine the $R$-stabilizers of the
action.  Observe that in this case, $W_R(e)=$ (in Osin's notation)
$B(R)=\{\gamma\in \Gamma\colon |\gamma|_{S\cup\sH}\le R\}.$

We extract the following results from the proof of \cite[Lemma
12]{Os}.

\begin{Proposition} For any $n\in\mathbb{N}$ define
$B(n)=\{\g\in \G\colon |\g|_{S\cup\sH}\le n\}.$ View
$B(n)\subset\Gamma$ with the metric $d_S.$ Fix $0\le\ell\le m.$
Then,
\begin{enumerate}
    \item the set $B(n-1)H_\ell$ can be written as
the disjoint union $B(n-1)H_\ell=\sqcup_{\g}\g H_\lambda,$ where
$g$ ranges over a certain (finite) subset of $B(n-1);$ and
    \item for any $r>0$ there is a subset $Y_r\subset \G$ so that
$Y_r=\bigcup_{x}B(n-1)x$ is a finite union and such that the sets
$\g H_\ell\setminus Y_r$ and $\g'H_\ell\setminus Y_r$ are
$r$-separated whenever $\g\neq \g'.$
\end{enumerate}
\end{Proposition}

Next, we need an analog of \cite[Lemma 12]{Os}:

\begin{Lemma} \label{OsinL12}  Suppose that there is a function
$f$ so that $\ad_{H_\ell}\le f$ for all subgroups $H_\ell.$  Then
for any $n\in\mathbb{N}$ we have $\ad_{B(n)}\preceq f.$
\end{Lemma}

\begin{proof}  This argument is essentially the one that appears
in \cite{Os}.

We proceed inductively. Observe that the fact that $|S|$ is finite
means that $B(1)= S\cup(\bigcup_{\ell=1}^mH_\ell)$ is a finite
union of sets with $\ad\preceq f.$ Thus, $\ad_{B(1)}\preceq f.$

For $n>1,$ we can write
$B(n)=(\bigcup_{\ell=1}^mB(n-1)H_\ell)\cup(\bigcup_{s\in
S}B(n-1)s).$ Since the set $B(n-1)s$ is quasi-isometric to
$B(n-1),$ we have $\ad_{B(n-1)s}\approx \ad_{B(n-1)}.$  Since
$|S|<\infty,$ the finite union theorem means that it remains only
to show that $\ad_{B(n-1)H_\ell}\preceq f$ for each $1\le\ell\le
m.$

By the previous proposition, $B(n-1)H_\ell=\sqcup_\g\g H\ell$ and
the growth of $\ad$ for $\{\g H_\ell\}$ is uniform in $\g$ since
these cosets are isometric.  By the remarks following Corollary
\ref{unioncor} we obtain $\ad_{B(n)}\preceq f$ as desired.
\end{proof}

We are now in a position to prove the main theorem from this
section.

\begin{Theorem} \label{relhyper} Let $\Gamma$ be a finitely generated group that is hyperbolic
relative to a finite family $\{H_1,\ldots,H_k\}$ of subgroups.
Suppose that there is some upper bound
$f:\mathbb{R}_+\to\mathbb{R}_+$ for $\ad_{H_i}$ for all $i.$ Then,
$\ad_\Gamma\preceq f.$
\end{Theorem}

\begin{proof} As observed above, $\Gamma$ acts by isometries on a
metric space with finite asymptotic dimension.  By Theorem
\ref{main} we need only check that for every $R>0$ there is some
$f:\mathbb{R}_+\to\mathbb{R}_+$ so that $\ad_{W_R}\preceq f.$ But,
this is exactly the content of Lemma \ref{OsinL12}.
\end{proof}

%
%
%
%
%

\section{Open Questions}

Dranishnikov \cite{Dr04} posed the following question, which could
be called the Milnor-type question for growth of $\ad:$

\begin{Question} Do there exist groups with intermediate dimension
growth?
\end{Question}

Also in \cite{Dr04}, Dranishnikov showed that restricted wreath
products of certain groups with finite asymptotic dimension have
polynomial dimension growth. That leads to a few natural
questions.

\begin{Question} Is the growth of the asymptotic dimension of a
restricted wreath product of groups with finite asymptotic
dimension at most polynomial?
\end{Question}

A positive answer would lead to:

\begin{Question} Is the class of groups with polynomial asymptotic
dimension growth closed under the formation of restricted wreath
products?
\end{Question}

A negative answer would resolve Question 1 and lead to:

\begin{Question} Is there a ``critical rate of growth" so that
groups whose dimension grows slower than this rate are exact and
those with faster growth are not? What about a critical rate for
coarse embeddings in Hilbert space?
\end{Question}

\begin{Question} Is it true that $\ad_G\approx\max\{\ad_H,\ad_K\}$
for an exact sequence $1\to K\to G\to H\to 1$ of finitely
generated groups?
\end{Question}


\providecommand{\bysame}{\leavevmode\hbox
to3em{\hrulefill}\thinspace}
\providecommand{\MR}{\relax\ifhmode\unskip\space\fi MR }
\providecommand{\MRhref}[2]{%
  \href{http://www.ams.org/mathscinet-getitem?mr=#1}{#2}
} \providecommand{\href}[2]{#2}

\end{document}